\documentclass[a4paper,12pt]{article}

\usepackage{amsbsy,eucal, amssymb,amsmath}
\usepackage{latexsym,amsfonts,amsthm}
\usepackage{graphicx}

\input xy
\xyoption{all}

\newcommand{\set}[1]{\{#1\}}
\newcommand{\Set}{{\sf Set}}
\newcommand{\Ord}{{\sf PreOrd}}
\newcommand{\Pos}{{\sf Pos}}
\newcommand{\Rel}{{\sf EndoRel}}
\newcommand{\Q}{{\sf Q}}

\newcommand{\pair}[1]{\langle #1 \rangle}
\theoremstyle{definition}

\newtheorem{definition}{\bf Definition}[subsection]
\newtheorem{lemma}{Lemma} [subsection]
\newtheorem{theorem}{Theorem} [subsection]

\newtheorem{exercise}{Exercise} [subsection]
\begin{document}
\title{Basic constructions in the categories of sets, sets with a binary relation on them, preorders and  posets}
\author{Ignacio Viglizzo\\
INMABB UNS CONICET, Departamento de Matemática\\
Bah\'ia Blanca, Argentina\\
{\tt viglizzo@gmail.com}}

\maketitle

\begin{abstract}
	The purpose of this note is to work out the details of the concrete incarnation of a few categorical constructions (products, coproducts, pullbacks, pushouts, equalizers, coequalizers and exponentials) in some useful and basic categories: the categories of sets, sets endowed with a binary relation, preorders and posets.
\end{abstract}

\section{Introduction}

This note came out of my need to work our some details I could not find in the textbooks on category theory. They are probably too trivial to be addressed there, but I still felt that they may prove useful for someone. I am not giving the definitions of categories, limits or functors for they can easily be found elsewhere. In what follows I have tried to be as clear as possible. One tool for doing so is the notation. Sometimes I am maybe too explicit, in order to remind of the nature of the objects we deal with. One example of this, is an attempt of distinguishing an element $a$ in a set $A$ from its inclusion $\pair{a,0}$ in a set $A+B$, but I think that at some point this fastidiousness pays off. 

Another thread through this note is to try to find an explanation why the proposed constructions turn out to be the ones described by the categorical definitions: why does a product have the number of elements it has and no other? Could a different binary relation fulfill the requirements to be part of the coproduct?

If you find these notes useful, have questions about them, or find errors to correct, please let me know.

We consider in this note the categories:
\begin{itemize}
	\item $\Set$, the category of sets and functions.
	\item  $\Rel$, the category of sets endowed with a binary relation on them (an \textit{endorelation}), with morphisms the set functions that preserve the binary relation. This is, $f:\pair{A,R}\to \pair{B,S}$ is a morphism in $\Rel$ if for every $a,a'\in A$, if $\pair{a,a'}\in R$ then $\pair{f(a),f(a')}\in S$.
	
	Note that while this category is named $\sf Rel$  in \cite{adamek06joy}. This category is also known as {\sf Bin}. Many other authors use $\sf Rel$ to denote the category of sets in which the morphisms are relations between sets. 
\item  $\Ord$%\footnote{This is the name given in the Wikipedia.}
, which has as objects pre-ordered sets, and monotonic functions as morphisms. This category is called $\sf Prost$ in \cite{adamek06joy}. 
\item $\sf Pos$ which has as objects posets, and as morphisms, monotonic functions among them. 
\end{itemize}

Just to be clear, a \textit{pre-order relation} over a set $A$ is a binary, reflexive and transitive relation. A preorder is a set with a pre-order relation over it, this is a pair $\pair{A, \preceq_A}$. Whenever there is no chance on confusion (and sometimes, even when there is), we will omit the subindex in the preorder relation, or simply  write that $A$ is a Preorder. A \textit{monotonic} (for preorders) function $f:\pair{A,\preceq_A}\to\pair{B,\preceq_B}$ is a set function $f:A\to B$ such that for every $a, b\in A$, if $a\preceq_Ab$ then $f(a)\preceq_Bf(b)$.

 Similarly, a \textit{poset} is a set $P$, endowed with a binary relation $\le_P$, which is an \textit{order} (also known as partial order): it is reflexive, transitive,  and antisymmetric. A function $f:\pair{P,\le_P}\to\pair{Q,\le_Q}$ is \textit{monotonic} if for every $x,y\in P$, $x\le_P y$ implies $f(x)\le_Q f(y)$.  As with preorders, we will often write the order relation as $\le$, leaving out the subindex that indicates which poset we are working on.

\section{$\Set$} 

\subsection{Products in $\Set$}

We will not repeat in this note the definitions of category or product in a category.  One of the basic properties of the products is that they are uniquely determined \textit{up to isomorphism}. Isomorphisms in $\Set$ are just bijective functions, so we might think that the important feature in a product of sets is the quantity of its elements. We usually realize the product of two sets $A$ and $B$ as their cartesian product:
\[A\times B=\set{\pair{a,b}: a\in A \text{ and } b\in B}.\]
So consider a simple example in which we have the finite sets $A=\set{a,b,c}$ and $B=\set{1,2}$. The six element set $A\times B=$ $\set{\pair{a,1},$ $\pair{a,2},$ $\pair{b,1},$ $\pair{b,2},\pair{c,1},\pair{c,2}}$, together with the implicit projection functions $\pi_A, \pi_B$ is \textit{a} product, but so is any other set with six elements, for example $A'=\set{\text{Black Widow, Captain America,Iron Man, Hawkeye, Hulk, Thor}}$, with adequate projection functions. 

What in the definition of products makes it so that it has exactly six elements? Suppose that a set $F$ with five elements is a product of $A$ and $B$ with projection functions $p_A:F\to A$ and $p_B:F\to B$ respectively. Now consider the set $A\times B$ with the usual projections $\pi_A$ and $\pi_B$. Since $F$ is a product, there must exist a function $u:A\times B\to F$ such that the following diagram commutes:
\[\xymatrix{
	&A\times B\ar[dl]_{\pi_A}\ar[dr]^{\pi_B}\ar@{-->}[d]^{u}& \\
	A&F\ar[l]^{p_A}\ar[r]_{p_B}&B
}\]
Since $A\times B$ has six elements and $F$ five, there must exist pairs $\pair{x,y}\neq\pair{z,t}\in A\times B$ such that $u(\pair{x,y})=u(\pair{z,t})$. We must then have that either $x\neq z$ or $y\neq t$, while on the other hand $x=\pi_A\pair{x,y}=p_Au\pair{x,y}=p_Au\pair{z,t}=\pi_A\pair{z,t}=z$ and similarly we can prove that $y=t$. 

Could a product  have \textit{more} than six elements? Assume now that a set $S$ with seven elements, with projections $p_A$ and $p_B$ is a product. We then have as above, a \textit{unique} function $u:A\times B\to S$. Since $S$ has seven elements there is at least two elements $s, s'$ with $p_A(s)=p_A(s')$ and $p_B(s)=p_B(s')$. The element $r=u\pair{p_A(s),p_B(s)}\in S$ is one of such elements,  and let $r'\neq r$ be another one. We can now define $u':A\times B\to S$ by 
 \[ u'(x)=\begin{cases}
 r'& \text{if } x=\pair{p_A(s),p_B(s)},\\
 u(x)& \text{otherwise}.
 \end{cases} \] In this way, we have a function $u'\neq u$ that makes the diagram commute, contradicting the uniqueness of $u$.

\subsection{Coproducts in $\Set$}

If $A$ and $B$ are sets, consider the \textit{disjoint union} of $A$ and $B$, which we will denote with $A+B$ and can be realized as $A\times\set{0}\cup B\times\set{1}$. Let $i_A:A\to A+B$ be the inclusion function $i_A(a)=\pair{a,0}$ for each $a\in A$, and similarly, $i_B:B\to A+B$ be the function $i_B(b)=\pair{b,1}$ for each $b\in B$.

We can take the same two sets from the example for products and ask again, why does $A+B$ have exactly five elements?  Suppose we had another set $C$ with functions $inc_A:A\to C$ and $inc_B:B\to C$ such that $C$ is a coproduct. Then we have in particular a unique function $[i_A,i_B]$:

\[\xymatrix{
	A\ar[r]^{inc_A}\ar[dr]_{i_A}	&C\ar@{-->}[d]^{[i_A,i_B]}	&B\ar[l]_{inc_B}\ar[dl]^{i_B}\\
	& A+ B					 	& \\
}
\]

If $C$ has more than five elements, then there are elements in $C$ that are not in the image of $A$ or $B$, so it doesn't matter where do we send them in $A+B$, so we can construct different functions making the diagram commute, and we lose the uniqueness of $[i_A,i_B]$. If $C$ has less than five elements, then we have that either for some  $a\in A$ and $b\in B$, $inc_A(a)=inc_B(b)$, or for some $a, a'\in A, a\neq a'$, $inc_A(a)=inc_A(a')$ (or we have a similar situation for $B$). In the first  case, 
 we must also have $\pair{a,0}=i_A(a)=[i_A,i_B]inc_A(a)=[i_A,i_B]inc_B(b)=i_B(b)=\pair{b,1}$, which is an impossibility.  In the second case we get $\pair{a,0}=\pair{a',0}$, which is an impossibility as well.\footnote{We have used here the fact that $i_A=\pair{\cdot,0}$ is an injective function. This is apparent from the definition we have given for the inclusion maps.  More generally, one can prove that the coproduct inclusions are monic morphisms if in the category the product is distributive with respect to the coproduct. See \cite{trimbleInjective}.}

\subsection{Pullbacks in $\sf Set$}

If $f:{A}\to {D}$ and $g:{B}\to {D}$ are functions with common codomain, let $P=\set{\pair{a,b}\in A\times B: f(a)=g(b)}$.  The set ${P}$, together with the projections $\pi_A$ and $\pi_B$, is a pullback in $\Set$. 

It is clear that the projections make the square commute:

\[\xymatrix{
	P\ar[r]^{\pi_A}\ar[d]_{\pi_B}	&A\ar[d]^{f}	\\
	B\ar[r]_{g}& D					 	
}
\]

Now let $X$ be a set  and consider functions  $h:{X}\to{A}$ and $j:X\to B$ such that 
\[\xymatrix{
	X\ar[rrd]^{h}\ar[rdd]_{j} &  &  \\
	&	&A\ar[d]^{f}	\\
	&	B\ar[r]_{g}& D					 	
}
\]
commutes. For each $x\in X$, since $f\circ h(x)=g\circ j(x)$, the pair $\pair{h(x),j(x)}$ is in $P$. Define then $u:X\to P$ by $u(x)=\pair{h(x),j(x)}$ for every $x\in X$.

\[\xymatrix{
	X\ar[rrd]^{h}\ar[rdd]_{j} 	\ar@{-->}[dr]^{u}&  &  \\
	&	P\ar[r]_{\pi_A}\ar[d]^{\pi_B}&A\ar[d]^{f}	\\
	&	B\ar[r]_{g}& D					 	
}
\]
For the uniqueness of $u$, consider a function $v$ that makes the diagram commute as well. Then $\pi_Av(x)=h(x)$ and $\pi_Bv(x)=j(x)$, so $v(x)=\pair{h(x),j(x)}=u(x)$.

\subsection{Pushouts in $\sf Set$}

Now let $C$ be a set, with morphisms (that is, functions) $f:C\to A$ and $g:C\to B$. We want to build the pushout of this diagram. To do this, first consider the coproduct set $A+B$.

On the set $A+B$ we now consider the least equivalence relation that contains  the pairs $\pair{f(c),g(c)}$ for all the elements $c\in C$.

An explicit way of characterizing the least equivalence relation $R^e$ containing a binary relation $R$ is as follows:

\begin{lemma} \label{ler}
	The least equivalence relation  containing a binary relation $R$ on a set $A$ consists of the pairs $\pair{a,b}\in A\times A$ such that there exist $n\ge 0$ and $a_0, a_1, \ldots a_n\in A$ such that $a=a_0, b=a_n$, and for each $i=1, \ldots , n$, either $\pair{a_{i-1}, a_{i}}\in R$ or $\pair{ a_{i}, a_{i-1}}\in R$. 
\end{lemma}
\begin{proof}
	Let $R^e$ be the relation defined in the statement. We need to prove that $R^e$ is an equivalence relation, that it contains $R$ and that it is the least such relation. 
	
	Taking $n=0$, we see that all the pairs $\pair{a,a}$ for $a\in A$ are in $R^e$, so $R^e$ is reflexive. 
	
	If $\pair{a,b}\in R^e$ then there exist $n\ge 0$ and $a_0, a_1, \ldots a_n\in A$ such that $a=a_0, b=a_n$, and for each $i=1, \ldots , n$, either $\pair{a_{i-1}, a_{i}}\in R$ or $\pair{ a_{i}, a_{i-1}}\in R$. Taking the sequence $b_0=a_n, b_1=a_{n-1},\ldots, b_n=a_0$ we readily see that $\pair{b,a}\in R^e$.
	
	If $\pair{a,b}\in R^e$  and $\pair{b,c}\in R^e$, then there exist $n\ge 0$ and $a_0, a_1, \ldots a_n\in A$ such that $a=a_0, b=a_n$, such that for each $i=1, \ldots , n$, either $\pair{a_{i-1}, a_{i}}\in R$ or $\pair{ a_{i}, a_{i-1}}\in R$ and  there exist $m\ge 0$ and $b_0, b_1, \ldots b_m\in A$ such that $b=b_0, c=b_m$ such that for each $j=1, \ldots , m$, either $\pair{b_{j-1}, b_{j}}\in R$ or $\pair{ b_{j}, b_{j-1}}\in R$. Taking the sequence $c_0=a, c_1=a_1,\ldots, c_n=a_n=b=b_0,c_{n+1}=b_1,\ldots, c_{n+m}=b_m=c$ we readily see that $\pair{a,c}\in R^e$.
	
	For each pair $\pair{a,b}\in R$ the sequence $a_0=a,a_1=b$ proves that $\pair{a,b}\in R^e$, so $R\subseteq R^e$.
	
	Finally, if $E$ is an equivalence relation that contains $R$ and $\pair{a,b}\in R^e$ then there exist $n\ge 0$ and $a_0, a_1, \ldots a_n\in A$ such that $a=a_0, b=a_n$, and for each $i=1, \ldots , n$, either $\pair{a_{i-1}, a_{i}}\in R$ or $\pair{ a_{i}, a_{i-1}}\in R$.  Since $E$ is symmetric and transitive, it turns out that all the pairs $\pair{a_{i-1}, a_{i}}$ are in $E$ and so is $\pair{a,b}$.
	
\end{proof}  

So our next step is to apply the definition of $R^e$ above to the case of the relation on $A+B$ consisting of the pairs $\pair{f(c),g(c)}$ for all the elements $c\in C$. It is the set $E$ of the pairs $\pair{x,y}\in (A+B)\times (A+B)$ such that there exist $n\ge 0$ and $x_0, x_1, \ldots x_n\in A+B$ such that $x=x_0, y=x_n$, and for each $i=1, \ldots , n$, either $\pair{x_{i-1}, x_{i}}=\pair{f(c),g(c)}$ for some $c\in C$ or $\pair{ x_{i-1}, x_{i}}=\pair{g(c),f(c)}$ for some $c\in C$. Let $Q$ be the quotient set $A+B$ over the equivalence relation $R^e$, and let $[\cdot]:A+B\to Q$ be the canonical quotient map.

\begin{lemma}\label{pushoutSet}
	The diagram  	
	\[\xymatrix{
		C\ar[r]^{f}\ar[d]_{g}	&A\ar[d]^{[\pair{\cdot,0}]}	\\
		B\ar[r]_{[\pair{\cdot,1}]}& Q					 	
	}
	\]
	
	is a pushout square.
\end{lemma}

\begin{proof} 
	
	We first observe that since for each $c\in C$, $f(c)$ and $g(c)$ are in the same equivalence class of $R^e$, $[\pair{f(c),0}]=[\pair{g(c),1}]$ so the square commutes.

	The elements of $Q$ are of the form $[\pair{f(c),0}]=[\pair{g(c),1}]$ for some $c\in C$ or they are singletons $[\pair{a,0}]=\set{\pair{a,0}}$, or $[\pair{b,1}]=\set{\pair{b,1}}$ if they come from elements in $A+B$ that are not in the image of $f$ or $g$.
	
	If there is a set $Y$ with functions $k:A\to Y$ and $l:B\to Y$ such that $k\circ f=l\circ g$, then we can define $v:Q\to Y$ by 
	
	\[ v(q)=\begin{cases}
	{k(f(c))}=l(g(c))& \text{if } q=[\pair{f(c),0}]=[\pair{g(c),1}] \text{ for some }c\in C,\\
	k(a)& \text{if } q=[\pair{a,0}] \text{ for some } a\notin f(C), \\
	l(b)& \text{if } q=[\pair{b,1}] \text{ for some } b\notin g(C).
	\end{cases} \]
	
	\[\xymatrix{
		C\ar[r]^{f}\ar[d]_{g}	&A\ar[d]_{[\pair{\cdot,0}]}\ar[ddr]^k	\\
		B\ar[r]^{[\pair{\cdot,1}]}\ar[drr]_l& Q\ar@{-->}[dr]^v				\\
		&&Y	 	
	}
	\]

	Here we need to check that $v$ is well-defined. If for example, $[\pair{f(c),0}]=[\pair{f(c'),0}]$, then there exist $n\ge 0$ and $x_0, x_1, \ldots x_n\in A+B$ such that $x_0=\pair{f(c),0}, x_n=\pair{f(c'),0}$, and for each $i=1, \ldots , n$, either $\pair{x_{i-1}, x_{i}}=\pair{\pair{f(d_i),0},\pair{g(d_i),1}}$ or $\pair{x_{i-1}, x_{i}}=\pair{\pair{g(d_i),1},\pair{f(d_i),0}}$ for some $d_i\in C$. Thus we get a sequence \[x_0=\pair{f(c),0}=\pair{f(d_0),0}.\]
	\[ x_1=\pair{g(d_0),1}=\pair{g(d_1),1},\]
	\[x_2=\pair{f(d_1),0}=\pair{f(d_2),0},\]
	\[\vdots\]
	\[x_n=\pair{f(d_{n-1}),0}=\pair{f(c'),0}.\]
	That is, for $i=1,\ldots,n$, 
	\[ x_i=\pair{g(d_{i-1}),1}=\pair{g(d_i),1}, \text{ if $i$ is odd, and }\]
	\[x_i=\pair{f(d_{i-1}),0}=\pair{f(d_i),0} \text{ if  $i$ is even. }\]
	Since the inclusions are injective, we get in particular that $g(d_{i-1})=g(d_i)$  if $i$ is odd and $f(d_{i-1})=f(d_i)$ if $i$ is even.
	
	But then  $k(f(c))=k(f(d_0))=l(g(d_0))=l(g(d_1))=k(f(d_1))=\ldots=k(f(c'))$.
	
	Notice that in the case above, $n$ must be even. If we had started assuming that $[\pair{f(c),0}]=[\pair{g(c'),1}]$, $n$ would have turned out to be odd.

	It is straightforward to check that $v\circ [\pair{\cdot,0}]=k$ and $v\circ [\pair{\cdot,1}]=l$. If $j:Q\to Y$ is such that $j\circ [\pair{\cdot,0}]=k$ and $j\circ [\pair{\cdot,1}]=l$, then $j$ must agree with the definition of $v$, so $v$ is the only function satisfying these two equations.
	
\end{proof}

\subsection{Equalizers and Coequalizers in $\sf Set$}

\begin{definition}
	Given two morphisms $f,g:X\to Y$, the equalizer is an object $E$ with a morphism $e:E\to X$ such that $f\circ e=g\circ e$ and if any morphism $m:Q\to X$ is such that $f\circ m=g\circ m$ then there exists a unique morphism $u:Q\to E$ such that $e\circ u=m$.
	
	\[\xymatrix{
		E\ar[r]^{e}&X \ar@<-.5ex>[r]_g \ar@<.5ex>[r]^f & Y	\\
		Q\ar@{-->}[u]_{u}\ar[ur]_m
	}
	\]
	
\end{definition}

In the case of sets, the equalizer of two functions $f,g:X\to Y$ can easily be described as the subset of $X$ in which the functions $f$ and $g$ agree. Let $E=\set{x\in X: f(x)=g(x)}$. Let $e:E\to X$ be the inclusion map. If there is a function $m:Q\to X$ such that $f\circ m=g\circ m$ then observe that for every $q\in Q$, $f(m(q))=g(m(q))$ so $m(q)\in E$, and we can put $u(q)=m(q)$. If $v:Q\to X$ is such that $e\circ v=m$, since $e$ is the identity for elements in $E$, we conclude that $v=m=u$.

\begin{definition}
	Given two morphisms $f,g:X\to Y$, the co-equalizer is an object $C$ with a morphism $c:Y\to C$ such that $c\circ f=c\circ g$ and if any morphism $m:Y\to P$ is such that $m\circ f=m\circ g$ then there exists a unique morphism $u:C\to P$ such that $u\circ c=m$.
	\[\xymatrix{
		X \ar@<-.5ex>[r]_g \ar@<.5ex>[r]^f & Y\ar[r]^c\ar[dr]_m&C\ar@{-->}[d]_{u}\\
		&&P
	}
	\] 	
\end{definition}

\begin{lemma}
	In the category of sets, consider the equivalence relation induced on $Y$ by the pairs $\pair{f(x),g(x)}$ for every $x\in X$. By Lemma \ref{ler} this relation consists of all pairs $\pair{y,y'}\in Y\times Y$ such that there exist $n\ge 0$ and $y_0, y_1, \ldots y_n\in Y$ such that $y=y_0, y'=y_n$, and for each $i=1, \ldots , n$, either $\pair{y_{i-1}, y_{i}}=\pair{f(x),g(x)}$ or $\pair{ y_{i-1}, y_{i}}=\pair{g(x),f(x)}$ for some $x\in X$. Taking $C$ to be the quotient of $Y$ by this relation and the canonical map from $Y$ to $C$ which we denote with $[\cdot]$, we obtain a co-equalizer of $f$ and $g$.
\end{lemma}

\begin{proof}
	It is immediate that for every $x\in X$, $[f(x)]=[g(x)]$. 
	
	Given $P$ and $m:Y\to P$ such that  $m\circ f=m\circ g$, we put for each $[y]\in C$, $u([y])=m(y)$. We need to check that $u:C\to P$ is a well defined function. If $[y]=[y']$ then there exist $n\ge 0$ and $y_0, y_1, \ldots y_n\in Y$ such that $y=y_0, y'=y_n$, and for each $i=1, \ldots , n$, either $\pair{y_{i-1}, y_{i}}=\pair{f(x),g(x)}$ or $\pair{ y_{i-1}, y_{i}}=\pair{g(x),f(x)}$ for some $x\in X$. For $n=0$, we have that $y=y'$, so $m(y)=m(y')$. If $n\ge 1$, then for each $i=1,\ldots,n$ we have that $[y_{i-1}]=[y_i]$. If the pair $\pair{ y_{i-1}, y_{i}}$ is $\pair{f(x),g(x)}$ or $\pair{g(x),f(x)}$ for some $x\in X$, then $m(f(x))=m(g(x))$,  so $m(y_{i-1})=m(y_i)$. It follows that $m(y)=m(y')$.
	
	If for some function $v:C\to P$, $v\circ[\cdot]=m$, then for each $[y]\in C$, $v([y])=m(y)$, so $v=u$.
\end{proof}

\subsection{Equalizers and Pullbacks}\label{sectequa}

In this section we work in arbitrary categories, and we see from general principles, that we have spent some unnecessary effort working out either the equalizer or pullback constructions in $\Set$. On the other hand, one of our goals is to have a concrete description of these constructions in particular categories.

There is a close relation between  equalizers  and pullbacks in categories in which binary products exist:

\begin{lemma}
	In a category with binary products and pullbacks, equalizers exist.
\end{lemma}

\begin{proof}
	Given $f,g:A\to B$, consider the diagram:
	
	\[\xymatrix{
		&B\ar[d]^{\Delta}	\\
		A\ar[r]_{\pair{f,g}}& B\times B					 	
	}
	\] where $\Delta=\pair{1_B,1_B}$. Take the pullback of this diagram. It will be an object $E$ and morphisms $e:E\to A$ and $h:E\to B$.
	The object $E$ and the morphism $e$ form a equalizer of $f$ and $g$.
	\[\xymatrix{
		E\ar[r]^h\ar[d]_e	&B\ar[d]^{\Delta}	\\
		A\ar[r]_{\pair{f,g}}& B\times B					 	
	}
	\]
	First we check that $f\circ e=g\circ e$. Let $\pi_1, \pi_2$ be the projections from $B\times B$ to $B$. Then we have $\pi_1\pair{f,g}e=fe=\pi_1\pair{1_b,1_b}h=h$. On the other hand, $\pi_2\pair{f,g}e=ge=\pi_2\pair{1_b,1_b}h=h$, so $fe=h=ge$.
	
	Now let $m:Q\to A$ be such that $fm=gm$. Then the square 
	\[\xymatrix{
		Q\ar[r]^{fm}\ar[d]_m	&B\ar[d]^{\Delta}	\\
		A\ar[r]_{\pair{f,g}}& B\times B					 	
	}
	\] commutes. By the definition of the pullback, there exists a unique $u:Q\to E$ such that $eu=m$, which is is the condition required to prove that $E$ and $e$ form an equalizer.
\end{proof}

\begin{lemma}
	In a category with binary products and equalizers, pullbacks exist.
\end{lemma}

\begin{proof}
	Consider the diagram:
	\[\xymatrix{
		&B\ar[d]^{g}	\\
		A\ar[r]_{f}&C					 	
	}
	\] We can construct the product $A\times B$ and then take the equalizer $E$ of the parallel arrows from $A\times B$ to $C$ $f\pi_A$ and $g\pi_B$.
	
	\[\xymatrix{
		E\ar[r]^{e}&A\times B \ar@<-.5ex>[r]_{g\pi_B} \ar@<.5ex>[r]^{f\pi_A} & C
	}
	\] We now check that the square 
	\[\xymatrix{
		E\ar[r]^{\pi_Be}\ar[d]_{\pi_Ae}&B\ar[d]^{g}	\\
		A\ar[r]_{f}&C					 	
	}\] is a pullback. It is clear that the square commutes because of the equalizer equation $f\pi_Ae=g\pi_Be$. 
	
	Now if $X$ is such that there are morphisms $k:X\to A$ and $l:X\to B$ such that $fk=gl$, then we have that $f\pi_A\pair{k,l}=g\pi_B\pair{k,l}$, so by the definition of equalizers, there exists a unique $u:X\to E$ such that $eu=\pair{k,l}$. 
	\[\xymatrix{
		E\ar[r]^{e}&A\times B \ar@<-.5ex>[r]_{g\pi_B} \ar@<.5ex>[r]^{f\pi_A} & C	\\
		X\ar@{-->}[u]_{u}\ar[ur]_{\pair{k,l}}
	}\ \ \ \ \ \ \xymatrix{
		X\ar[rrd]^{k}\ar[rdd]_{l} 	\ar@{-->}[dr]^{u}&  &  \\
		&	E\ar[r]_{\pi_Ae}\ar[d]^{\pi_Be}&A\ar[d]^{f}	\\
		&	B\ar[r]_{g}& C					 	
	}
	\]
	Therefore, $\pi_A eu=\pi_A\pair{k,l}=k$ and $\pi_B eu=\pi_B\pair{k,l}=l$.
\end{proof}

Furthermore, there is a connection between pullbacks and products. A terminal object \textbf{1} in a category can be seen as an empty product. If there is a terminal object in a category, then the pullback of any diagram of the form 
\[\xymatrix{
	&B\ar[d]^{!_B}	\\
	A\ar[r]_{!_A}&\textbf{1}					 	
}
\] is a product of $A$ and $B$.

Dually, if co-products exist in a category, pushouts can be constructed from coequalizers and coequalizers can be constructed from pushouts.

\subsection{Exponentials}

We just briefly review the usual construction of exponentials in $\Set$. Given sets $A$ and $B$, let $B^A=\set{f:A\to B}$ that is, the set of all the functions from $A$ to $B$. 

The set $B^A$ can also be denoted as $A\Rightarrow B$, a notation that is perhaps more suggestive of the set of functions going from $A$ to $B$, and also agrees with the intuitionistic implication in Heyting algebras of sets.

The function $eval: B^A\times A\to B$ is defined by $eval(f,a)=f(a)$ for every $a\in A$ and $f:A\to B$. Now we can check the defining universal property of exponentials, namely that for every function $f:X\times A\to B$, there is a unique function $\tilde{f}:X\to B^A$ such that the following diagram commutes:

\[\xymatrix{
B^A&B^A\times A\ar[r]^{eval}&B\\
X\ar[u]^{\tilde{f}}&X\times A\ar[u]^{\tilde{f}\times 1_A}\ar[ur]_f&
}\]

For each $x\in X$, $\tilde{f}(x)$ is a function from $A$ to $B$ defined by $\tilde{f}(x)(a)=f(x,a)$ for all $x\in X$  and $a\in A$. Furthermore, for any function $g:X\to B^A$ such that $eval\circ g\times 1_A=f$ and all $x\in X, a\in A$  we have that $g(x)(a)=eval(g(x),a)=f(x,a)=\tilde{f}(x)(a)$.

We are interested as before in seeing why the defined set of all functions from $A$ to $B$ fits exactly with the definition of the exponential object. It is helpful to think in terms of the \textit{adjunction} of the functors $\cdot\times A$ and $\cdot^A$. The adjunction is an isomorphism (so, in $\Set$, a bijection) between the sets of all the functions from $X\times A$  $B$ and the set of all functions from $X$ to $B^A$. This isomorphism is usually written as:
\[\frac{X\times A\to B}{X\to B^A}\]

Now consider the special case in which $X$ is a singleton set $1=\set{*}$. Now we have:
\[\frac{1\times A\to B}{1\to B^A}\]
This is, we have a bijection between the elements in $B^A$ and the functions from $A$ to $B$.

\section{$\Rel$}

\subsection{Products in $\Rel$}

If $\pair{A,R}$ and $\pair{B,S}$ are sets with binary relations, we take the set $A\times B$ (the cartesian product), together with the relation $T$ defined by:
\[\pair{a,b}T\pair{a',b'}\text{ iff } aRa'\text{ and }bSb'\]
Clearly $T$ is a binary relation on $A\times B$. The projections $\pi_A:A\times B\to A$ and $\pi_B:A\times B\to B$ preserve the relations: if $\pair{a,b}T\pair{a',b'}$ then by the definition of $T$, $aRa'$ and $bSb'$.

Given a set $C$ with a relation $U$ on $C$, and morphisms $f:C\to A$, and $g:C\to B$, define $\pair{f,g}:C\to A\times B$ by:
\[\pair{f,g}(c)=\pair{f(c),g(c)}\]
for every $c\in C$.  Is $\pair{f,g}$ a morphism in $\Rel$? If $cUc'$ for $c,c'\in C$, then $f(c)Rf(c')$ and $g(c)Sg(c')$. Therefore, 
$\pair{f(c),g(c)}T\pair {f(c'),g(c')}$ this is $\pair{f,g}(c)T\pair{f,g}(c')$.

 From the definitions,  $\pi_A\circ \pair{f,g}=f$ and $\pi_B\circ \pair{f,g}=g$. Furthermore, if $k:C\to A\times B$ is a morphism such that $\pi_A\circ k=f$ and $\pi_B\circ k=g$, then for every $c\in C$, $k(c)=\pair{f(c),g(c)}$ so $k=\pair{f,g}$.

\[\xymatrix{
	&\pair{C,U}\ar[dl]_f\ar[dr]^g\ar@{-->}[d]^{\pair{f,g}}& \\
	\pair{A,R}&\pair{A\times B,T}\ar[l]^{\pi_A}\ar[r]_{\pi_B}&\pair{B,S}
}\]

The same reasoning is valid for building all \textit{small} products, that is, products over arbitrary sets.

Could a different relation $T'$ on $A\times B$ also be a product in $\Rel$? 
All the pairs in $T$ must be in $T'$: if $a,a'\in A$ and $b,b'\in B$ are such that  $\pair{a,a'}\in R$ and $\pair{b,b'}\in S$,  consider a set $C=\set{c,c'}$ with the relation $U=\set{\pair{c,c'}}$. We set $f:C\to A$ by $f(c)=a$ and $f(c')=a'$. In the same manner, $g:C\to B$ is given by $g(c)=b$ and $g(c')=b'$. Then the function $\pair{f,g}:C\to A\times B$ must preserve the relation $U$, so $\pair{f,g}(c)=\pair{a,b}$ should be related by $T'$ with $\pair{f,g}(c')=\pair{a',b'}$.

 If there was an excess of pairs in $T'$, there would exist  a pair $\pair{\pair{a,b},\pair{a',b'}}$ in $T'$, but not in $T$. Then we would have that $\pair{a,a'}\notin R$ or $\pair{b,b'}\notin S$. Say that $\pair{a,a'}\notin R$, then the projection $\pi_A$ would not preserve the relation $T'$. This is, assuming that the relations $R$ and $S$ are not empty! What happens when $S$ is empty?
 
 We may call the relation $T$ defined above $R\times S$.

\subsection{Coproducts in $\Rel$}

If $\pair{A,R}$ and $\pair{B,S}$ are sets with binary relations, we take the coproduct of their underlying sets  $A+B$  as we did in $\Set$.  We endow this set with the relation $V$:
\[\pair{x,y}V\pair{x',y'}\text{ iff } xRx' \text{ and } y=y'=0, \text{ or }xSx' \text{ and } y=y'=1.\]

Alternatively, if we identify the elements $a\in A$ with the pairs $\pair{a,0}$ and the elements $b\in B$ with the pairs $\pair{b,1}$, we can write:
\[xVy\text{ iff } (x, y\in A \land  xRy) \lor  (x, y \in B \land  xSy).\]

 Both $i_A$ and $i_B$ are readily seen to preserve the binary relation: if $aRa'$ for $a, a'\in A$ then $\pair{a,0}V\pair{a',0}$. We may also indicate the relation $V$ by $R+S$.

Given morphisms $f:\pair{A,R}\to \pair{D,W}$ and $g:\pair{B, S}\to \pair{D,W}$, we can use the function   $[f,g]:A+B\to D$ defined for the category $\Set$.

\begin{exercise}
	Check that $[f,g]$ is a morphism in $\Rel$.
\end{exercise}

We now argue that the given relation $R+S$ is the only one that make the set $A+B$ the coproduct in $\Rel$. Suppose $V'$ is a relation on $A+B$ such that $\pair{A+B,V'}$ is a coproduct, with inclusions $inc_A$ and $inc_b$. All the pairs from $R+S$ must be in $V'$: Take $A\neq\emptyset$ and $a,a'\in A$ such that $\pair{a,a'}\in R$. Since $i_A$ is a morphism in $\Rel$, $\pair{i_A(a),i_A(a')}\in V'$.  
  
  Now assume that for some $a,a'\in A$ we have that $\pair{\pair{a,0},\pair{a',0}}\in V$ but $\pair{a,a'}\notin R$. Then, the following diagram must commute:
\[\xymatrix{
	\pair{A,R}\ar[r]^{inc_A}\ar[dr]_{i_A}	&\pair{A+B,V'}\ar@{-->}[d]^{[i_A,i_B]}	&\pair{R,S}\ar[l]_{inc_B}\ar[dl]^{i_B}\\
	& \pair{A+ B,R+S}					 	& \\
}
\] Then, we have that $\pair{[i_A,i_B]\pair{a,0},[i_A,i_B]\pair{a',0}}\in R+S$, that is $\pair{i_A(a),i_A,(a')}\in R+S$, so by the definition of $R+S$, we get that $\pair{a,a'}\in R$, a contradiction.

\subsection{Pullbacks in $\Rel$}

If $f:\pair{A,R}\to \pair{D,W}$ and $g:\pair{B,S}\to \pair{D,W}$ are morphisms in $\Rel$ with common codomain, let $P=\set{\pair{a,b}\in A\times B: f(a)=g(b)}$, and define the relation $T$ by 
\[\pair{a,b}T\pair{a',b'}\text{ iff } aRa'\text{ and }bSb'.\]
This is the same relation as in the construction of the product. We claim that $\pair{P,T}$, together with the projections $\pi_A$ and $\pi_B$, is a pullback in $\Rel$. 

It is clear that the projections preserve the binary relations, and  the square 
\[\xymatrix{
	\pair{P,T}\ar[r]^{\pi_A}\ar[d]_{\pi_B}	&\pair{A,R}\ar[d]^{f}	\\
	\pair{B,S}\ar[r]_{g}& \pair{D,W}					 	
}
\] 
commutes.

Now let $\pair{X,V}$ be a set with a binary relation and  morphisms $h:\pair{X,V}\to\pair{A,R}$ and $j:\pair{X,V}\to \pair{B,S}$ such that 
\[\xymatrix{
	\pair{X,V}\ar[rrd]^{h}\ar[rdd]_{j} &  &  \\
	&	&\pair{A,R}\ar[d]^{f}	\\
	&	\pair{B,S}\ar[r]_{g}& \pair{D,W}					 	
}
\]
commutes. As we did in $\Set$, define then $u:X\to P$ by $u(x)=\pair{h(x),j(x)}$ for every $x\in X$. To check that $u$ is preserves the relation, assume $xVy$ for some $x, y \in X$. Since $h$ and $j$ preserve the relation, $h(x)Rh(y)$ and $j(x)Sj(y)$, so $\pair{h(x),j(x)}T\pair{h(y),j(y)}$, that is,  $u(x)Tu(y)$.

\[\xymatrix{
	\pair{X,V}\ar[rrd]^{h}\ar[rdd]_{j} 	\ar@{-->}[dr]^{u}&  &  \\
	&	\pair{P,T}\ar[r]_{\pi_A}\ar[d]^{\pi_B}&\pair{A,R}\ar[d]^{f}	\\
	&	\pair{B,S}\ar[r]_{g}& \pair{D,W}					 	
}
\]
The uniqueness of $u$ follows easily as in the product.

\subsection{Pushouts and Coequalizers in $\Rel$} \label{pushoutRel}

In order to construct pushouts in $\Rel$, we define the quotient by an equivalence relation $E$ of an endorelation $R$ on a set $A$ by:
\[R/E=\set{\pair{[a], [b]}:\pair{a,b}\in R}\]
This means that a pair $\pair{[a], [b]}$ is in $R/E$ if and only if there exist $a'\in A$ and $b'\in B$ such that $a'\in [a], b'\in[b]$ and $\pair{a',b'}\in R$. Notice that this is the \textit{least} endorelation on $A/E$ that makes the function $[\cdot]$ a morphism in $\Rel$: with this definition, we ensure that every time that  we have that $\pair{a,b}\in R$,  $\pair{[a], [b]}\in R/E$, and at the same time, no extra pair of elements appears in $R/E$.

To check that the construction of the pushout we gave for sets works here, we need to see that the function $v$ defined in Lemma \ref{pushoutSet} preserves the endorelation. For this, assume that we have a pair $\pair{[a],[b]]}\in R+S/E$. Then, there exist $a',b'$ such that $a'\in[a]$, $b'\in[b]$, and $\pair{a',b'}\in R+S$. Assuming $a'=\pair{c,1},b'=\pair{d,1}$ for some $c,d\in B$ such that $\pair{c,d}\in S$, it follows that $\pair{l(c),l(d)}\in U$, and this means that $\pair{v([a]),v([b])}\in U$, since $v([a])=v([a'])=v([\pair{c,1}])$ and $v([b])=v([b'])=v([\pair{d,1}])$.  Needless to say, the situation is similar if $a'$ and $b'$ come from elements in $A$.

	\[\xymatrix{
	\pair{C, U}\ar[rr]^{f}\ar[d]_{g}&	&\pair{A,R}\ar[d]_{[\pair{\cdot,0}]}\ar[ddrr]^k	\\
	\pair{B,S}\ar[rr]^{[\pair{\cdot,1}]}\ar[drrrr]_l&& \pair{A+B/E,R+S/E}\ar@{-->}[drr]^v				\\
	&&&&\pair{Y,U}	 	
}
\]

The coequalizers can be constructed from co-products and pushouts as shown in section \ref{sectequa}.

\subsection{Terminal object}

There is a terminal object in the category $\Rel$: take a singleton $\set{*}$ with the identity relation $\set{\pair{*,*}}$. Then, for every $\pair{A,R}$ the map $!_A$ that sends every element in $A$ to $*$ is a morphism in $\Rel$.

\subsection{Exponentials in $\Rel$}

We have seen that the exponential in $\Set$ of sets $A$ and $B$ can be constructed the set of all functions from  $A$ to $B$. In $\Rel$, both $A$ and $B$ are endowed with a binary relation, and we need to put a binary relation on this set accordingly. We will first give the relation, then we check that $eval$ and $\tilde{f}$ are morphisms in $\Rel$ for every morphism $f:X\times A\to B$, and finally we will justify why the relation defined is the only one that could work.

We write with $\pair{B,S}^{\pair{A,R}}$ or simply $B^A$  the set of all functions from $A$ to $B$. Given $\alpha, \beta\in \pair{B,S}^{\pair{A,R}}$, we define:
\[\pair{\alpha,\beta}\in W \text{ iff }\pair{a,a'}\in R\text{ implies } \pair{\alpha(a),\beta(a')}\in S.\]

$eval: \pair{B^A,W}\times\pair{A,R}\to\pair{B,S}$ is a morphism in $\Rel$: Assume $\pair{\pair{\alpha,a},\pair{\beta,a'}}\in W\times R$, so $\pair{\alpha,\beta}\in W$ and $\pair{a,a'}\in R$. Then, by the definition of $W$,  $\pair{\alpha(a),\beta(a')}\in S$,  that is, $\pair{eval(\pair{\alpha,a}),eval(\pair{\beta,a'})}$, just as we needed.

Now assume that $f:\pair{X,V}\times \pair{A,R}\to\pair{B,S}$ is a morphism, and $\tilde{f}:X\to B^A$ is defined as before by $\tilde{f}(x)(a)=f(x,a)$ for every $x\in X, a\in A$. Let $\pair{x,x'}\in V$. Then to prove that $\pair{\tilde{f}(x),\tilde{f}(x')}\in W$, we need to take a pair $\pair{a,a'}\in R$ and prove that $\pair{\tilde{f}(x)(a),\tilde{f}(x')(a')}\in S$. But this is straightforward because from the hypothesis $\pair{\pair{x,a},\pair{x',a'}}\in V\times R$, so $\pair{f({x,a}),f({x',a'})}\in S$ which is $\pair{\tilde{f}(x)(a),\tilde{f}(x')(a')}\in S$.

To analize the relation $W$ we will use again the adjunction, now with $2=\set{0,1}$ with the relation $V=\set{\pair{0,1}}$. The product $\pair{2,V}\times\pair{A,R}$ can be interpreted as a pair of copies of $A$: $\set{0}\times A$ and $\set{1}\times A$, and the relation $V\times R$, for each pair $\pair{a,a'}\in R$, has exactly one pair joining $\pair{0,a}$ with $\pair{1,a'}$. 

At the same time, a morphism $p$ from $\pair{2,V}$ to any $\pair{C,T}$ has the effect of picking one element of the relation $T$: $p(0)$ is the first component, and $p(1)$ the second, while $p$ is a morphism in $\Rel$ only if $\pair{p(0),p(1)}\in T$.
\[\frac{\pair{2\times A,V\times R}\to \pair{B,S}}{\pair{2,V}\to \pair{B^A,W}}\]
So, there is a bijection between the morphisms from $\pair{2,V}$ to $\pair{B^A,W}$, and morphisms from $\pair{2\times A,V\times R}$ to $\pair{B,S}$. The former correspond to pairs in $W$, while there are as many in the latter as morphisms from $\pair{ A, R}$ to $\pair{B,S}$.

\subsection{$\Rel$ is unbalanced}

Consider a set with two different endorelations on it. To be more concrete, take the set $A=\set{a,b}$ and the objects $\pair{A,\emptyset}$ and $\pair{A,\set{\pair{a,b}}}$. Then the identity is a $\Rel$-morphism from the first to the second, which is monic and epic but has no inverse. Categories in which this can happen are called \textit{unbalanced}. 

\section{Preordered sets}

$\Ord$ is a full subcategory of $\Rel$. This means that all the morphisms in $\Rel$ that go from a preordered set to another are also monotonic functions and therefore morphisms in $\Rel$. Because of this, we don't need to check that the functions used with the constructions in $\Rel$ are monotonic.

\subsection{Products in $\Ord$}

If $\pair{A,\preceq_A}$ and $\pair{B,\preceq_B}$ are preoders, we take the set $A\times B$ (the cartesian product), together with the relation $\preceq_{A\times B}$ defined by:
\[\pair{a,b}\preceq_{A\times B}\pair{a',b'}\text{ iff } a\preceq_Aa'\text{ and }b\preceq_Bb'\]

This is a particular case of the one for $\Rel$. We should check that $\preceq_{A\times B}$ is a  preorder relation, but this is straightforward. It also follows from the work in $\Rel$ that the projections $\pi_A:A\times B\to A$ and $\pi_B:A\times B\to B$ are monotonic, and that given a preordered set $C$, and morphisms $f:C\to A$, and $g:C\to B$, the function $\pair{f,g}:C\to A\times B$  defined by $\pair{f,g}(c)=\pair{f(c),g(c)}$ for every $c\in C$ is the unique morphism such that $\pi_A\pair{f,g}=f$ and $\pi_B\pair{f,g}=g$ 

\subsection{Coproducts in $\Ord$}
If $\pair{A,\preceq_A}$ and $\pair{B,\preceq_B}$ are preorders, consider the construction $\pair{A+B, \preceq_A+\preceq_B}$. We may write $\preceq_A+\preceq_B$ as $\preceq_{A+B}$.
\[\pair{x,y}\preceq_{A+B} \pair{x',y'}\text{ iff } x\preceq_Ax' \text{ and } y=y'=0, \text{ or }x\preceq_Bx' \text{ and } y=y'=1.\]
The functions $i_A:A\to A+B$  and $i_B:B\to A+B$ work as in $\Rel$ and checking the universality of the construction requires no new considerations.

\subsection{Pullbacks in $\Ord$}

If $f:A\to D$ and $g:B\to D$ are morphisms in $\sf Ord$ with common codomain, the pullback of this diagram can be constructed as in $\Rel$, after checking that the restriction of the preorder relation on the product to the subset $P=\set{\pair{a,b}\in A\times B: f(a)=g(b)}$ is still a preorder.
\[\xymatrix{
	P\ar[r]^{\pi_A}\ar[d]_{\pi_B}	&A\ar[d]^{f}	\\
	B\ar[r]_{g}& D					 	
}
\]

\subsection{Pushouts in $\sf Ord$}

To construct a pushout or co-equalizer in $\Ord$, we face a problem: the quotient of a preordered set by an arbitrary equivalence relation may not be preordered. In particular, transitivity fails. Consider the example on the following figure:

	\includegraphics{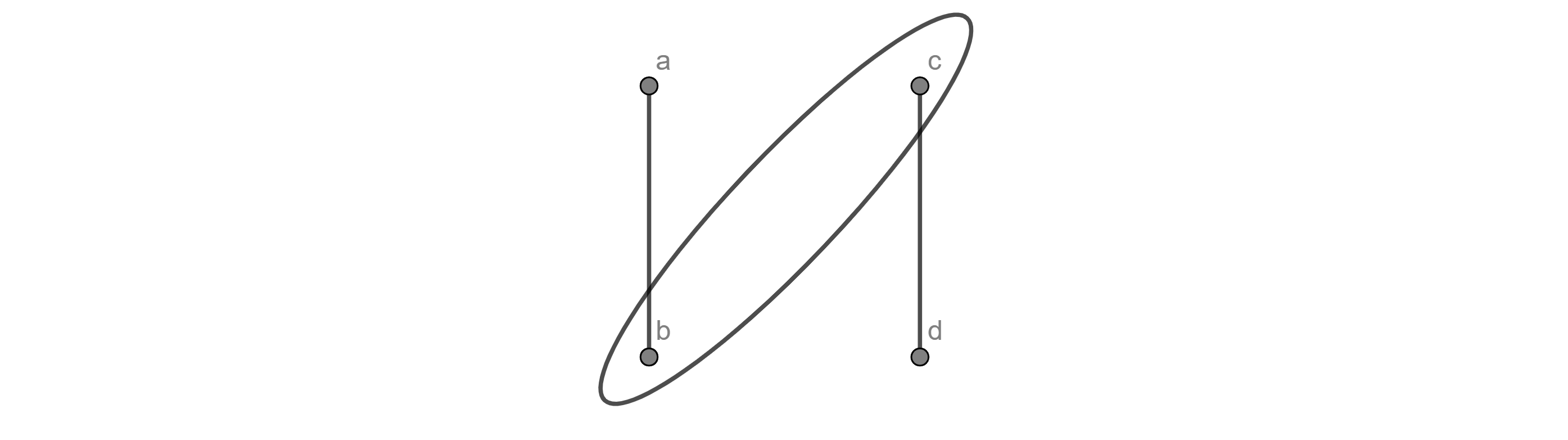}

We have a preordered set and an equivalence relation with classes $\set{a}, \set{b,c}$, and $\set{d}$. The quotient set would have $[b]\preceq[a]$ and $[d]\preceq[c]=[b]$ but no relation between $[a]$ and $[d]$. 

When we constructed a quotient in $\Rel$, we endowed the quotient set with the least relation that made the canonical map a morphism. Now we must find the least relation on the quotient set that makes it a preorder, while still keeps the quotient map a morphism. For this, we will use the  \textit{transitive closure} of a binary relation.

\begin{lemma}
	Given a binary relation $R$ on a set $A$,  the least transitive relation $R^t$ containing $R$ is the set of pairs $\pair{a,b}\in A\times A$ such that there exist $a_1, a_2, \ldots a_n$ such that $a=a_1, b=a_n$, and for each $i=1, \ldots , n-1$,  $\pair{a_i, a_{i+1}}\in R$. 
\end{lemma} 

The proof of the lemma above is similar to the one of Lemma \ref{ler}.

Now let $\pair{A,\preceq_A}, \pair{B,\preceq_B}, \pair{C,\preceq_C}$ be preordered sets, with morphisms $f:C\to A$ and $g:C\to B$. We want to build the pushout of this diagram. To do this, first consider the coproduct preorder $\pair{A+B, \preceq_A+\preceq_B}$. As we did in $\Rel$, we  consider  on $A+B$ the equivalence relation $R^e$ where 
\[R=\set{\pair{f(c),g(c)}:c\in C}  \]
 Then we take the quotient set with the corresponding relation $\pair{A+B/_{R^e},\preceq_{A+B}/_{R^e} }$ and finally we take the transitive closure of this relation, to obtain the preordered set $\pair{A+B/_{R^e},(\preceq_{A+B}/_{R^e})^t }$.

\begin{exercise}
	If $R$ is a  reflexive relation on a set $A$, then its transitive closure, $\preceq$ is a preorder relation on $A$.
\end{exercise}

We will write  $\preceq$ for $(\preceq_{A+B}/_{R^e})^t$. Since $\preceq_{A+B}$ is reflexive,  so is $\preceq_{A+B}/_{R^e}$ and therefore $\preceq$ is a preorder on $A+B/_{R^e}$. We will also write $(A+B/_{R^e})^t$ for the pair $\pair{A+B/_{R^e},\preceq}$. 
%Notice that the identity map from $A+B/_{R^e}$ to  $(A+B/_{R^e})^t$ is a monotone function.

\begin{lemma}
$(A+B/_{R^e})^t$, together with the functions $[\cdot]i_A$ and $[\cdot]i_B$, is a pushout in $\Ord$. 
\end{lemma} 

\begin{proof}
\[\xymatrix{
	C\ar[r]^{f}\ar[d]_{g}	&A\ar[d]_{[\cdot]i_A}\ar[ddr]^k	\\
	B\ar[r]^{[\cdot]i_B\ \ }\ar[drr]_l& (A+B/_{R^e})^t\ar@{-->}[dr]^v\\
	&&Y					 	
}
\]
We use the same functions that we used for $\Rel$ in the subsection \ref{pushoutRel}. The point that remains to be checked is that the function $v$ defined is a morphism when considered from the preordered set $(A+B/_{R^e})^t$ to $Y$. 	

If we consider $[x]\preceq[y]$, then there exist $[x_1]=[x], [x_2], \ldots,[x_n]\in A+B/_{R^e}$ such that for each $i=2,\ldots,n$, $\pair{[x_{i-1}],[x_i]}\in \preceq_{A+B}/_{R^e}$. This means that for each $i=2,\ldots,n$, there exist $y_{i-1}\in[x_{i-1}], y_i\in[x_i]$ such that $y_{i-1}, y_i\in A$ and $y_{i-1}\preceq_A y_i$ or $y_{i-1}, y_i\in B$ and $y_{i-1}\preceq_B y_i$. In the first case, we get that $v([x_{i-1}])=k(y_{i-1})\preceq_Yk(y_i)=v([x_i])$ and in the second, $v([x_{i-1}])=l(y_{i-1})\preceq_Yl(y_i)=v([x_i])$. Using the transitivity of $\preceq_Y$, we obtain $v([x])\preceq_Yv([y])$. Notice that each $[x_i]$ may have representatives in the set $A$ and $B$ simultaneously.
\end{proof}

%\subsection{Further results about $\Ord$}
%
%$\Ord$ is a 2-category (see wikipedia). It is a cartesian closed category (nLab). It is an exponential ideal in the cartesian closed category $\sf Cat$.

\section{The category  $\sf Pos$ of posets}

\subsection{$\Pos$ is a reflective subcategory of $\Ord$}

Given a preorder relation $\preceq$ on a set $A$, one can define an equivalence relation $\approx$ on $A$ by $a\approx b$ if and only if $a\preceq b$ and $b\preceq a$. If we denote with $[a]$ the equivalence class, then the relation on $A/\approx$ defined by $[a]\le[b]$ if and only if there exist $a'\in [a]$ and $b'\in[b]$ such that $a'\preceq b'$ is a partial order, so $\pair{A/\approx,\le}$ is a poset.

First note that $\Pos$ is a full subcategory of $\Ord$. We will prove that $\Pos$ is a also a \textit{reflective} subcategory of $\Ord$. The definition is taken from \cite{adamek06joy}:

\begin{definition}
	Let {\sf A} 	be a subcategory 	of {\sf B}, and 	let $B$ 	be a {\sf B}-object. 
	\begin{enumerate}
		\item  	An {\sf A}-\textit{reflection} (or {\sf A}-\textit{reflection 	arrow}) for $B$ 	is a {\sf B}-morphism	$r:B\to A$ from 	$B$ to an {\sf A}-object $A$ with the following universal property:
		
	for any {\sf B}-morphism $f:B\to A'$ from 	$B$ into some {\sf A}-object $A'$, there exists a unique 	{\sf A}-morphism 	$f': A \to  A'$ such that the triangle
	\[\xymatrix{
		B\ar[r]^{r}\ar[dr]_{f}	&A\ar@{-->}[d]^{f'}	\\
		& A'
	}
	\]
	commutes. By an ``abuse of language'' an object $A$ is called an 	{\sf A}-reflection 	for $B$ provided that there exists an {\sf	A}-reflection
$r:B\to  A$ for $B$ with 	codomain 	$A$.
	\item 	{\sf A} 	is 	called 	a \textit{reflective 	subcategory} 	of {\sf B} 	provided 	that each 	{\sf B}-object 	has an {\sf A}-reflection.
	
\end{enumerate}
\end{definition} 

\begin{theorem}
	$\Pos$ is a reflective subcategory of $\Ord$.
\end{theorem}
\begin{proof}
	Given a preordered set $\pair{P,\preceq_P}$, define the relation $\equiv$ by 
	\[p\equiv p' \text{ iff } p\preceq_P p'\text{ and }p'\preceq_Pp.\]
	$\equiv$ is an equivalence relation and $\pair{P/_\equiv,\preceq/_\equiv}$ is a poset. The canonical map $[\cdot]:P\to P/_\equiv$ is a {\sf Pos}-reflection.
	
	Let $\pair{A,\leq_A}$ be a poset and $f:P\to A$ a morphism in {\sf Ord}. Define $f':P/_\equiv\to A$ by $f'([p])=f(p)$ for each $[p]\in P/_\equiv$.
	
		\[\xymatrix{
		P\ar[r]^{[\cdot]}\ar[dr]_{f}	&P/_\equiv\ar@{-->}[d]^{f'}	\\
		& A
	}
	\]
	
	$f'$ is well-defined: If $[p]=[p']$ then $p\preceq_P p'\text{ and }p'\preceq_Pp$. Then $f(p)\le_A f(p')$ and $f(p')\le_A f(p)$, so $f(p)=f(p')$.

	$f'$ is a morphism in $\Ord$: if $[p]\preceq/_\equiv[q]$, then there exist $p'\in[p]$ and $q'\in[q]$ such that $p'\preceq_Pq'$ so $f'([p])=f(p')\le_Af(q')=f'([q])$.
	
	Finally $f'$ is unique: if $g:P/\equiv\to A$ is such that $g\circ[\cdot]=f$, then for each $p\in P$, $g([p])=f(p)=f'([p])$.
\end{proof}

\begin{exercise} Prove that in the Theorem above, $\equiv$ is an equivalence relation and $\pair{P/_\equiv,\preceq/_\equiv}$ is a poset.
\end{exercise}

Some general results from reflective subcategories apply then, so we have:

\begin{lemma}
	There is a unique functor ${\sf Q:Ord\to Pos}$ such that for every preordered set $\pair{A,\preceq_A}$, $\Q(A)=\pair{A/_\equiv,\preceq_A/_\equiv}$ and for each morphism $f:\pair{A,\preceq_A}\to\pair{B, \preceq_B}$ in $\Ord$, the diagram 
		\[\xymatrix{
		A\ar[d]^{f}\ar[r]^{q_A}	&\Q(A)\ar[d]^{\Q(f)}	\\
		B\ar[r]_{q_B}& \Q(B)
	}
	\] commutes, where $q_A:A\to\Q(A)$ is the canonical quotient map (or the $\Pos$-reflection arrow).
\end{lemma}

The proof can be seen as a particular case of 4.22 in \cite{adamek06joy}, or proved as an exercise.

\subsection{Products, coproducts and pullbacks}
The constructions of products, coproducts and pullbacks in $\Ord$, when applied to posets, yield posets, so these can be constructed in $\Pos$ just as in $\Ord$ and in $\Rel$, and their universal properties are preserved.

\subsection{Coequalizers in $\sf Pos$}

\begin{lemma} \cite{adamek06joy} 7.69, p. 119. In $\Pos$ coequalizers are formed by first forming the coequalizers in $\Ord$ and then taking the $\Pos$-reflection of the resulting preordered set.
\end{lemma}

\begin{proof}
	Let $f,g:X\to Y$ be $\Pos$-morphisms. We can consider them in the bigger category $\Ord$ and their coequalizer there, $c:Y\to C$. 
	
	Now let $Z$ be a poset with a $\Pos$-morphism $m:Y\to Z$ such that $mf=mg$, and let $q:C\to Q=\Q(C)$ be the $\Pos$ reflection of $C$.
	\[\xymatrix{
	X \ar@<-.5ex>[r]_g \ar@<.5ex>[r]^f & Y\ar[r]^c\ar[drr]_m&C\ar@{-->}[dr]^{u}\ar[r]^q&Q \ar@{-->}[d]^{v}\\
	&&&Z
}
\] 		
By the definition of coequalizer, there exists a unique $\Ord$-morphism $u:C\to Z$ such that $uc=m$, and by the definition of $\Pos$-reflection, there exists a unique $\Pos$-morphism $v:Q\to Z$ such that $vq=u$. It follows that $vqc=uc=m$, and if $g$ is a $\Pos$ morphism such that $gqc=m$, then since $vqc=m$, from the definition of reflection it turns out that $v=g$.
\end{proof}

\subsection*{Acknowledgment}
Thanks to Fernando Tohm\'e for the discussions and reading this manuscript.

\bibliographystyle{alpha}

\begin{thebibliography}{AHS06}
	
	\bibitem[AHS06]{adamek06joy}
	Ji\v{r}\'{i} Ad\'{a}mek, Horst Herrlich, and George~E. Strecker.
	\newblock Abstract and concrete categories: the joy of cats.
	\newblock {\em Repr. Theory Appl. Categ.}, (17):1--507, 2006.
	\newblock Reprint of the 1990 original [Wiley, New York; MR1051419].
	
	\bibitem[Tri15]{trimbleInjective}
	Todd Trimble.
	\newblock Distributivity implies monicity of coproduct inclusions, 2015.
	\newblock
	\url{https://ncatlab.org/toddtrimble/published/distributivity+implies+monicity+of+coproduct+inclusions}.
	
\end{thebibliography}

\end{document}